\documentclass{article}
\usepackage{epsf}
\usepackage[dvips]{graphics}

\newfont{\mm}{msbm10}

\begin{document}

\begin{titlepage}

\noindent
	{\LARGE {\bf Interval straight line fitting}}

\bigskip\noindent
{\large MAREK W. GUTOWSKI}
\par\noindent
{\em Institute of Physics, Polish Academy of Sciences, Warszawa, Poland\/}
\par\noindent
{\em e-mail: gutow@ifpan.edu.pl}

\bigskip\noindent
{\bf Abstract.} I consider the task of so called {\em curve fitting\/}
or {\em experimental data fitting\/} commonly encountered in various
branches of scientific research.  Unlike the traditional approach I do
not try to minimize any functional based on available experimental
information, instead the minimization problem is replaced with
constraint satisfaction procedure, which produces the interval hull of
solutions of desired type. The method, called {\em box slicing
algorithm\/}, is described in details. The results obtained this way
need not to be labeled with {\em confidence level\/} of any kind, they
are simply certain (guaranteed).  Additionally, the memory requirements
for the presented method are very conservative.  The approach is
directly applicable to other experimental data processing problems like
outliers detection and finding the straight line, which is tangent to
the experimental curve.

\bigskip\noindent
{\bf keywords:} experimental data fitting, curve fitting, interval
enclosures, experimental uncertainty, linear regression, constraint
satisfaction, box slicing, outlier identification, asymptotes

\vfil
{\small
\noindent
Address for correspondence:
\begin{quote}
Marek W. Gutowski\\
Institute of Physics, Polish Academy of Sciences\\
Al. Lotnik\'ow 32/46\\
02--668 Warszawa\\
Poland
\end{quote}
\noindent
e-mail:  gutow@ifpan.edu.pl
}
\end{titlepage}

\date{
}



\section{Stating the problem}
In many branches of science the so called problem of experimental data
fitting is often encountered.  Its short description is following:
\begin{quote}
Given the data set, i.e. the set of experimental observations, called
also measurements, $\left\{ \left(x_{j},\
y_{j}\right)\right\}_{j=1}^{n}$, and the model, find appropriate
parameters of this model, which adequately describe the data.
\end{quote}
In the following I will assume that:
\begin{itemize}
\item
the values of both coordinates for each measurement may be uncertain,
i.e. for each $x_{j}$ (respectively $y_{j}$) we know the interval,
$\left[\underline{x}_{j},\ \overline{x}_{j}\right]= {\mathbf x}_{j}$
(resp. $\left[\underline{y}_{j},\ \overline{y}_{j}\right]= {\mathbf
y}_{j}$), containing $x_{j}$ (resp. $y_{j}$) and guaranteed to contain
the true, unknown value of the measured or controlled physical
quantity,
\item
we are searching for parameters of a linear model relating $x$ with $y$:
$$
	y=ax+b
$$
\end{itemize}
The further considerations are directly applicable to many other
two-pa\-ra\-me\-ter models.  The extension for even more complicated
cases is also straightforward.  I choose the linear model since it is
very important, widely used and, at the same time, probably the
simplest one.

\medskip\noindent
Shortly: the goal is to find the bounds for two parameters, called
hereafter $a$ and $b$, which describe reliably the experimental data. 
There are many procedures to solve such a~problem, all depending on the
exact meaning what is the {\em best\/} solution, the LSQ (least squares
method) and LAD (least absolute deviations) being among the most
popular.  Yet, even the interval counterparts of those methods do not
deliver the results expected by experimentalists.  Often encountered is
the ``cluster problem'' (see \cite{Du1} and \cite{Du2}), which makes
the precise location of a global extremum difficult.  Additionally, as
a consequence of clustering, the enclosures for physically interesting
parameters are usually very pessimistic, up to the point of complete
unusability.  So why bother at all with one more interval method?

\section{Deficiencies of existing methods}

Most popular and commonly used fitting methods are nowadays the ones
based on probabilistic grounds.  This is because the results of
measurements are treated as random variables.  There is nothing wrong
with such an assumption, however further treatment of the experimental
data is most often than not based on, rarely explicitly stated,
additional assumptions concerning the distributions of measured
quantities.  It is assumed, and almost never checked, that the
distributions are normal (Gaussian).  However, contrary to common
belief, they usually aren't normal.  Today the vast majority of
measurements is performed with digital measuring devices, so even if
the investigated phenomenon is normally distributed, then the set of
its measurements, consisting of discrete values only, cannot be
normally distributed.

\medskip\noindent
There is also one more hypothesis being used, namely that the
uncertainties ({\em errors\/}) are {\em small\/}.  This is never
checked, and indeed cannot be checked, since there is no possibility to
influence the uncertainties of measurements once they had been
performed.

\medskip\noindent
Finally, all such methods, explicitly or implicitly,  make use of the
Central Limit Theorem, without ever bothering, that the conclusions
drawn on those grounds are only asymptotically valid, in the limit of
infinite number of measurements.

\medskip\noindent
The estimates of interesting parameters, obtained this way, are given
as a pair of numbers meaning: either the most probable value and its
standard deviation (again silently assuming the normal distribution!)
or --- less often --- as the confidence interval.  The choice of the so
called confidence level, which is then a third number, remains
arbitrary.  Needless to say, that the confidence level is only very
loosely related, if at all, to the performed measurements.

\section{Interval point of view}

In this paper we are going to find the tight and guaranteed bounds for
both parameters $a$ and $b$.  According to this aim, we will search for
the intervals,
${\mathbf a} =
\left[\underline{a},\ \overline{a}\right]$ and ${\mathbf b} =
\left[\underline{b},\ \overline{b}\right]$, containing {\em with
certainty\/} the true values of $a$ and $b$ respectively.  This is
equivalent with finding the solutions of the following system of
equations:
\begin{equation}\label{eq_system}
\left\{
\begin{array}{lcr}
	{\mathbf a}{\mathbf x}_{1} + {\mathbf b} & = & {\mathbf y}_{1} \nonumber\\
	\vdots                            & \vdots   & \vdots \nonumber\\
	{\mathbf a}{\mathbf x}_{n} + {\mathbf b} & = & {\mathbf y}_{n}
\end{array}
\right.
\end{equation}
It is a system of linear equations with $n>2$ equations and only $2$
unknowns. Since the number of available data exceeds the number of
unknowns, then the system (\ref{eq_system}) is overdetermined and
therefore generally has no solutions in the usual sense.  Nevertheless,
we will find such intervals ${\mathbf a}$ and ${\mathbf b}$, that the
equations (\ref{eq_system}) and the experimental data will be in some
sense consistent.  According to Shary \cite{Shary}, there are many ways
of saying what kinds of solutions are of interest to us. To count them
all let us rewrite the system (\ref{eq_system}) in standard matrix
form:
\begin{equation}\label{matrix_form}
\left(
\begin{array}{cc}
{\mathbf x}_{1} & 1\\
{\mathbf x}_{2} & 1\\
\vdots          & \vdots\\
{\mathbf x}_{n} & 1
\end{array}
\right)
\left(
\begin{array}{c}
{\mathbf a}\\
{\mathbf b}
\end{array}
\right) = \left(
\begin{array}{c}
{\mathbf y}_{1}\\
{\mathbf y}_{2}\\
\vdots\\
{\mathbf y}_{n}
\end{array}
\right)
\end{equation}
To define a particular set of solutions, we have to assign one of two
available quantifiers, $\forall$ or $\exists$, to each entry of the
matrix of coefficients of equations and to each component of the right
hand side vector (for details see \cite{Shary} and the following
discussion).  This makes $2^{2n}\times 2^{n}=2^{3n}$ possible
assignments corresponding to this many various solutions sets. 
Fortunately, this rather huge number can be substantially reduced. 
First, the second column of the matrix of coefficients consists of
simple real numbers only, not intervals.  There is no point in
assigning any quantifier to the numbers originating from the degenerate
interval $\left[1,\ 1\right]$, since there is exactly only one number
contained in it. Secondly, we do not want to distinguish any particular
measurement from among others.  This means, that the quantifiers have
to be assigned in a~special way: all ${\mathbf x}$'s should be
connected with {\em the same\/} quantifier.  The same can be said about
${\mathbf y}$'s.  Taking all this into account we arrive with only $4$
kinds of solutions of the system (\ref{eq_system}).  They are
following:
\begin{description}
\item [solutions in the usual sense:]
\begin{equation}
	\left\{\left(a,\ b\right):\quad \forall_{k=1\, \ldots\,
	n}\; \forall_{x\in{\mathbf x}_{k}}\; \forall_{y\in{\mathbf
	y}_{k}}\quad ax+b=y\right\}
\end{equation}
This set is unbounded when there is only one experimental point, exact
or not.  When the number of measurements is equal $2$ and the
measurements are exact (and different, i.e. ${\mathbf x}_{1}\cap
{\mathbf x}_{2} = \emptyset$), then it reduces to the single point. For
two inexact measurements it is bounded, when the measurements are
disjoint, and unbounded otherwise.  It is usually empty, when we have
$3$ or more measurements, since it is impossible to draw a straight
line connecting three or more arbitrarily chosen points, each belonging
to its own rectangle ${\mathbf x}_{k}\times {\mathbf y}_{k}$, even if
each of those rectangles consists of only one point.  The solutions in
usual sense are then of no interest for the experimentalists except,
perhaps, the case $n=2$, which is described later.  For now let's
assume that
$n>2$.
\item [united solutions:]
\begin{equation}
	\left\{\left(a,\ b\right):\quad \forall_{k=1\, \ldots\,
	n}\; \exists_{x\in{\mathbf x}_{k}}\; \exists_{y\in{\mathbf
	y}_{k}}\quad ax+b=y\right\}
\end{equation}
This set, if non-empty, consists by definition of all straight lines
having at least one common point with each of the rectangles ${\mathbf
x}_{1}\times {\mathbf y}_{1}$, ${\mathbf x}_{2}\times {\mathbf y}_{2}$,
\ldots\ ${\mathbf x}_{n}\times {\mathbf y}_{n}$.  In set theory
language the above may be expressed as
\begin{equation}\label{united}
	\forall_{k=1\, \ldots\, n}\;\left({\mathbf a}{\mathbf x}_{k} +
	{\mathbf b}\right) \cap {\mathbf y}_{k} \ne \emptyset
\end{equation}
Of course, not every pair $\left(a,\ b\right)\in\left({\mathbf a},\
{\mathbf b}\right)$ is the member of the set of solutions.  Indeed,
whenever we say ``solutions'', then we really mean ``interval hull of
the solution set''.  The same comment applies to the two remaining
cases described below.
\item [controllable solutions:]
\begin{equation}
	\left\{\left(a,\ b\right):\quad \forall_{k=1\, \ldots\,
	n}\; \exists_{x\in{\mathbf x}_{k}}\; \forall_{y\in{\mathbf
	y}_{k}}\quad ax+b=y\right\}
\end{equation}
Treating intervals as sets, we can write the relation satisfied by
interval hull of solutions of this kind as
\begin{equation}\label{control}
	\forall_{k=1\, \ldots\, n}\; {\mathbf a}{\mathbf x}_{k} +
	{\mathbf b} \supseteq {\mathbf y}_{k}
\end{equation}
\begin{figure}[h]
 \epsfysize=0.56\hsize
 \epsfxsize=0.76\hsize
 \centerline{\epsfbox{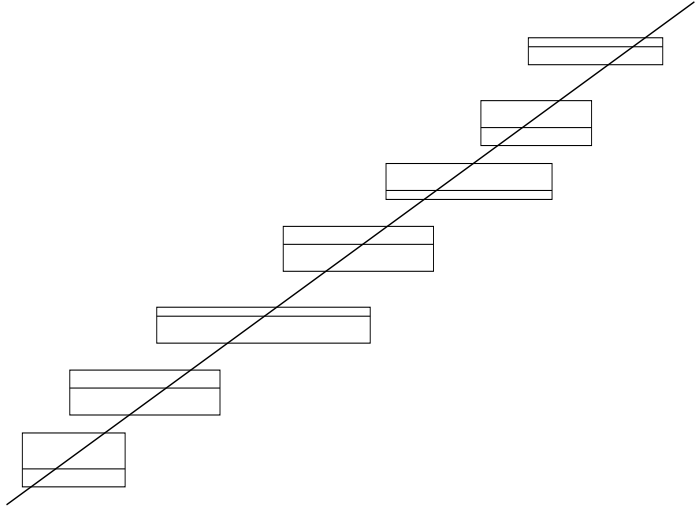}}
\caption{\sl Geometric interpretation of the {\bf controllable
solution}. The straight line shown belongs to the set of controllable
solutions, since it passes through all, arbitrarily chosen, short
horizontal sections inside each experimental rectangle.  It does {\bf
not} cross any of the vertical edges of any rectangle.
}
\label{con_solutions_fig}
\end{figure}

Finally, we should consider
\item [tolerable solutions:]
\begin{equation}
	\left\{\left(a,\ b\right):\quad \forall_{k=1\, \ldots\,
	n}\; \forall_{x\in{\mathbf x}_{k}}\; \exists_{y\in{\mathbf
	y}_{k}}\quad ax+b=y\right\}
\end{equation}
which are contained in the box $\left({\mathbf a},\ {\mathbf
b}\right)$, for which the relation
\begin{equation}\label{tolerance}
	\forall_{k=1\, \ldots\, n}\; {\mathbf a}{\mathbf x}_{k} +
	{\mathbf b} \subseteq {\mathbf y}_{k}
\end{equation}
holds.
\begin{figure}[h]
 \epsfysize=0.56\hsize
 \epsfxsize=0.76\hsize
 \centerline{\epsfbox{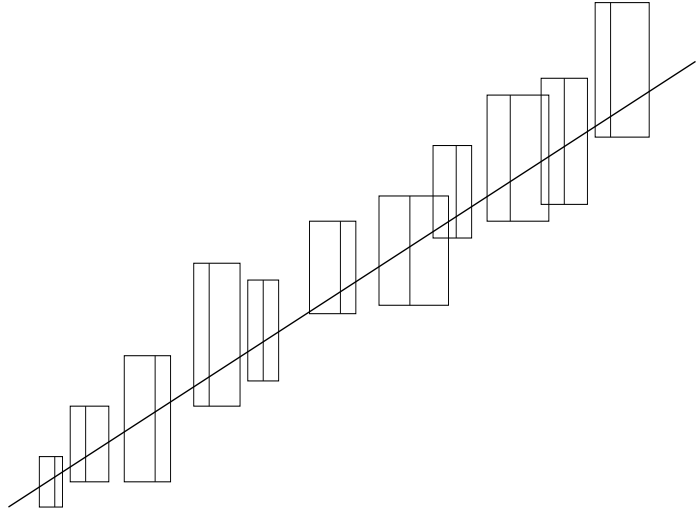}}
\caption{\sl Geometric interpretation of the {\bf tolerable solution}.
The straight line shown belongs to the set of tolerable solutions,
since it passes through all, arbitrarily chosen, short vertical
sections inside each experimental rectangle. It does {\bf not} cross
any horizontal edge of any rectangle.
}
\label{tol_solutions_fig}
\end{figure}
\end{description}
All the above relations should hold for every measurement; the names of
various sets of solutions are used after Shary \cite{Shary}.

\medskip\noindent
In the strictly mathematical sense we have covered all kinds of
solutions.  Experimentalists may be interested in one more type of
``solutions'', which will be called here {\bf crude solutions}.  These
are defined in almost the same way as the united solutions, except that
the relations (\ref{united}) are required to hold for majority of
measurements, not necessarily for all of them.  This way every united
solution is, of course, a crude solution as well.  Such ``solutions''
might be useful, when analyzing data containing {\em outliers\/}.  Let
us refrain from further discussion of crude solutions now, leaving it
to the later part of this paper. Instead consider still another set,
defined as the smallest one satisfying the conditions
\begin{equation}\label{trash}
	\forall_{k=1\, \ldots\, n}\; \forall_{x\in {\mathbf x}_{k}}\;
	{\mathbf a}x + {\mathbf b} \supseteq {\mathbf y}_{k}
\end{equation}
which are very similar to those for the interval hull of controllable
solutions.  Simply stating, the set just defined has the intuitively
simple property, that the graph of the expression
${\mathbf y}={\mathbf a}x + {\mathbf b}$ covers all experimental
uncertainty rectangles in the $xy$ plane.  I shall show now, that the
set (\ref{trash}) is unbounded in {\mm R}$^{2}$:

\medskip\noindent
Suppose that we have fixed parameter ${\mathbf a}=\left[a,\
a\right]\in\ ${\mm IR} as a thin (degenerate) interval being equal to
the arbitrary real number.  It is easy to adjust the parameter
${\mathbf b} = \left[\underline{b},\ \overline{b}\right]\in\ ${\mm IR}
such that (\ref{trash}) holds; it is enough to put
\begin{equation}\label{low_b}
	\underline{b}\ \le\ \min_{k}\ \biggl(\underline{\mathbf y}_{k} -
	\underline{a{\mathbf x}}_{k}\biggr)
\end{equation}
and
\begin{equation}\label{high_b}
	\overline{b}\ \ge\ \max_{k}\ \biggl(\overline{\mathbf y}_{k} -
	\overline{a{\mathbf x}}_{k}\biggr)
\end{equation}
The smallest ${\mathbf b}$ corresponds, of course, to the situation,
when we have equalities in (\ref{low_b}) and (\ref{high_b}).  But,
since $a$ was arbitrary, then the entire set is unbounded in {\mm
R}$^{2}$, and so must be its convex (interval) hull.  Thus we have
shown, that the smallest set of type (\ref{trash}), in the sense that
it is contained in any other having required properties, does not
exist.

\medskip\noindent
We will not explore the idea of finding the box with minimal volume,
since this might lead to unplausible or even unphysical results. 
Consider the poor quality experimental data, (for sake of simplicity we
assume that both variables are dimensionless, so the inequality below
makes sense) for which
\begin{equation}
	{\rm width}\left(\underline{\bigcup}\ {\mathbf x}_{k}\right)
	< {\rm width}\left(\underline{\bigcup}\ {\mathbf
	y}_{k}\right)
\end{equation}
and $\underline{\cup}$ denotes the convex hull.  For such data it may
easily happen, that the ``fitted'' line will be perpendicular to the
expected direction in the $xy$ plane!

\medskip\noindent
In conclusion, we will not discuss further the solutions of type
(\ref{trash}).  There is no hope, that they would constitute good
starting point for subsequent refinement to other types of solutions,
more tight.  It may even happen, that the randomly selected set of
``solutions'' with property (\ref{trash}) has no common elements with
the united solutions set.

\section{Properties and usability of remaining types of solutions}

First note that the following inclusions are always true:
\begin{equation}
	{\bf tolerable\ \ solutions} \subseteq {\bf united\ \
	solutions}
\end{equation}
and
\begin{equation}
	{\bf controllable\ \ solutions} \subseteq {\bf united\ \
	solutions}
\end{equation}
This means, that the set of united solutions is the ``largest'' one among
all of the true solutions (we exclude the crude solutions, since they
are not true in the sense, that they do not hold for every experimental
point).

\begin{figure}[h]
 \epsfysize=0.56\hsize
 \epsfxsize=0.76\hsize
 \centerline{\epsfbox{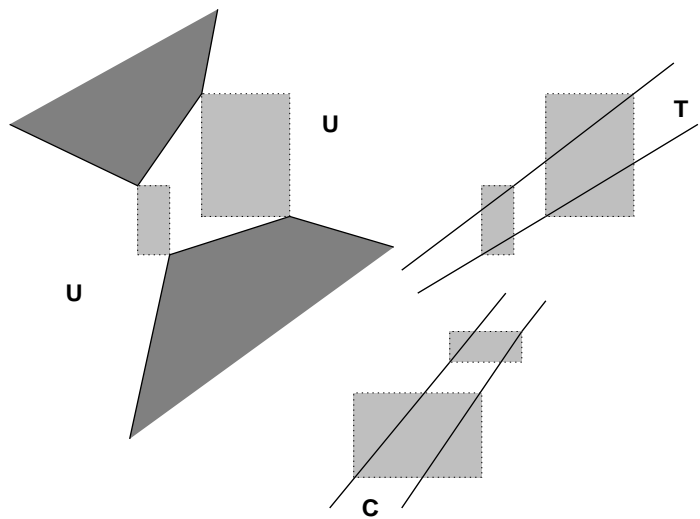}}
\caption{\sl United (U) vs. tolerable (T) and controllable (C) solutions
for two experimental points.  Uncertain data are shown as light grey
rectangles. Solid lines bound regions, marked with appropriate letters,
where the straight lines belonging to each solution set can be drawn. 
For united solutions only the dark shadowed region is inaccessible.
}
\label{solutions_fig}
\end{figure}
\medskip\noindent
The natural question arises, which kind of solutions (and why) is the
one, which should be used by experimentalists?  The quick answer is
{\em they all deserve our attention\/}. Looking at
Fig.\ref{solutions_fig}, we can see that the sets of tolerable and/or
controllable solutions deliver potentially better, i.e. more accurate
(tight) and thus tempting estimates of unknown parameters.  Yet, for
the experimentalist they are unreliable.  Suppose that after some time
the new experimental result is available, obtained with much better
accuracy.  It may place itself within the one of the already known
uncertainty rectangles but {\em outside} the domain marked with letter
{\sf T} (or {\sf C}) in Fig.\ref{solutions_fig}.  It is obvious, that
in such circumstances the tolerable (or controllable) solution set will
be empty. No definite conclusion can be drawn whether or not our
knowledge had really increased as a result of this new measurement. 
Driven by the widely accepted paradigm that increase of knowledge is
nothing else as decrease of ignorance, (and this should never increase
after new, correctly performed measurement) we would rather prefer the
united solutions' set.  Indeed, the choice between various kinds of
solutions is rather limited, since usually either tolerable or
controllable solutions will be available, if any, but not both.  It is
quite obvious, that significant progress, i.e. better bounds for
unknown parameters, can only be achieved when new measurements, of
similar accuracy, are taken off the already investigated range of
controlled parameter ${\mathbf x}$.  Repeating measurements, within
already explored range, is unlikely to improve significantly the bounds
of searched parameters, unless they are definitely more accurate.

\medskip\noindent
Summarizing the above considerations I propose to interpret various
sets of solutions in a following uniform way:
\begin{itemize} 
\item if {\bf united solutions exist}, then the available data are in
agreement with the model in use; there is no apparent contradiction
between theory and data.  The wording ``fair'', ``satisfactory'',
``good'' or ``excellent agreement'' is a matter of taste rather than
anything else.  Using of such phrases may be only justified by
comparison with similar results, especially concerning the widths of
intervals ${\mathbf a}$ and ${\mathbf b}$, obtained by different
method(s) or by other authors.
\item if {\bf no united solution exist} (two other sets of solutions
are then empty too) then either of the following happened:
\begin{itemize}
\item
one or more points are unreliable, i.e. their error bounds are
underestimated, perhaps even all of them, or
\item
one or more points are {\em outliers\/}.  This may be the result of
malfunctioning apparatus, errors in data transmission, or simply human
mistake when writing down the instrument's readings or typing them into
computer.
\item the linear model is not applicable to the phenomenon under
investigation.
\end{itemize}

The latter may happen to be true quite easily in physical sciences,
where the approximate, linearized models are used frequently.  They are
usable only as long as newer and better, more accurate, results become
available.  Then it might be the time to revise, correct or even reject
the current theory.

We can also see the united solution set from other perspective:  if it
is empty, then we have a~{\bf proof}, that our model is not adequate
for the observed data.  Inadequacy of the model may be bad news, but,
on the other hand, the proof of this fact is much more valuable than
the result of any statistical test.  For this to be true, we have to be
sure, that the uncertainties of all our data were estimated correctly,
never underestimated.
\end{itemize}

\medskip\noindent
Few words should be said concerning two other types of solutions. The
existence of either tolerable or controllable solutions, in addition to
the united set of solutions, has no peculiar meaning, at least for
experimentalists.  The lack of tolerable solutions may suggest that the
uncertainties connected with variable $x$, which is usually under
control, are either overestimated or should be reduced in further
experiments, perhaps with better apparatus.  The tolerable solutions
may only exist, when those uncertainties are small enough.  When all
$x$'s are exact, then tolerable solutions are simply united solutions. 
Similarly, the controllable solutions exist only when the tolerances
for all $y$'s are tight enough.  Again, when those tolerances are equal
to zero, then the controllable solutions are simply united solutions.
So, the existence or lack of existence of tolerable or controllable
solutions may be regarded only as a hint of how well the measurements
were performed.  As the final outcome of an experiment we should,
however, use {\em only\/} the intervals ${\mathbf a}$ and ${\mathbf b}$
derived from a united solution set, no matter that the solutions of
other kind usually look better.

\section{The algorithm}

In this section I present the general strategy of finding the interval
enclosure of every kind of solutions of (\ref{eq_system}).  It may be
regarded as a functional counterpart of the interval functions {\sf
ZERO1} and {\sf ZERO2} defined by van~Emden in \cite{Emden}. The method
is called {\em box slicing\/} or {\em box peeling algorithm\/}.

\medskip\noindent
The first step is to convert (\ref{eq_system}) into equivalent set of
conditions, usually in a form of inequalities.  We also need to
determine an initial box $V$, containing all the potential solutions of
desired type.  Let's defer the discussion on the exact forms of those
conditions and the choices of initial box for later. Our goal is to
find the smallest interval box $\left({\mathbf a},\ {\mathbf
b}\right)\in$ {\mm IR}$^{2}$, containing all pairs $(a,\ b)$, for which
appropriate conditions are satisfied. The general outline of the
algorithm is following:

\medskip
\fbox{
\begin{minipage}[h]{0.9\hsize}
\begin{center}
{\sf Input}
\end{center}
Initial box $V \in$ {\mm IR}$^{2}$ containing all solutions.

\medskip
\begin{center}
{\sf Algorithm}
\end{center}
\noindent
For each unknown interval in turn do the following:
\begin{itemize}
\item try to slice the box $V$ from the left; replace $V$ with new
box, if~slicing was successful
\item try to slice the box $V$ from the right; replace $V$ with new
box, if~slicing was successful
\end{itemize}
If any slicing, in any unknown,  was successful, then repeat the
procedure.

\medskip
\begin{center}
{\sf Output}
\end{center}
\noindent
Tight interval hull $\left({\mathbf a},\ {\mathbf b}\right)$ for
parameters $a$ and $b$.

~
\end{minipage}
}

\subsection{What is slicing?}
Suppose that the current box is $V=\left({\mathbf p}_{1},\ {\mathbf
p}_{2},\ \ldots\ {\mathbf p}_{r}\right)\in$ {\mm IR}$^{r}$ and we are
currently working with parameter ${\mathbf p}_{k}$. Slicing from the
left is described as a~sequence of steps:
\begin{enumerate}
\item $\xi \leftarrow\ 1 $
\item\label{restart} $\xi \leftarrow\ \xi/2$
\item divide $V$ into two parts, by cutting it with the plane
$p_{k}=p_{\xi}$, where $p_{\xi} = \underline{\mathbf p}_{k} +
\xi\left(\overline{\mathbf p}_{k} - \underline{\mathbf p}_{k}\right)$. 
Call the newly created subboxes {\bf slice\/} ($p_{k} \le p_{\xi}$)
and {\bf rest\/} ($p_{k} \ge p_{\xi}$).
\item probing \cite{Emden} {\bf slice} means determining, whether or
not the subbox {\bf slice} fails the considered system of inequalities.

If {\bf slice} fails the system of inequalities then
\begin{itemize}
\item $V \leftarrow\ $ {\bf rest} (discard {\bf slice})
\item finish slicing from the left with parameter ${\mathbf p}_{k}$;
exit with flag {\em success\/}
\end{itemize}
else (probe thinner slice)
\begin{itemize}
\item if termination criteria are not met then goto \ref{restart} else
finish slicing from the left with parameter ${\mathbf p}_{k}$; exit
with flag {\em no success\/}
\end{itemize}
\end{enumerate}

\medskip\noindent
Slicing from the right is similar, except that at the beginning $\xi$
is set to $0$ and the later updates have the form: $\xi \leftarrow
(1+\xi)/2$.

\medskip\noindent
Few comments are in order.  As seen from the description, at the
beginning the algorithm probes large chunks of the initial box $V$. 
Indeed, the first slice has the same volume as the remaining part of
box $V$, while the subsequent slices are smaller and smaller.  We stop
slicing in the direction of the current parameter at first success and
then immediately switch to the next parameter, according to van~Emden's
suggestions \cite{Emden}.

\medskip\noindent
What is the termination criterion?  The most obvious should be the one
based on width of a slice.  Unsuccessful slicing should stop, at the
latest, when the width of a slice, in direction of currently processed
parameter, becomes small, comparable with the machine accuracy.  One
might think, that we should terminate slicing in a given direction even
earlier, at some predefined threshold $\varepsilon$: slicing ceases,
when $\xi \le \varepsilon$ (slicing from the left) or $1-\xi \le
\varepsilon$ (slicing from the right), where $\varepsilon$ is small,
arbitrarily chosen, positive number, usually in range
$10^{-6}$---$10^{-3}$.  This choice, however, does not guarantee
obtaining the tightest possible bounds for searched parameters. 
Nevertheless, it may be practical in terms of of CPU time, and quite
sufficient when processing experimental data.

\medskip\noindent
Concluding, we may estimate the temporal complexity of a single cycle
of slicing all unknowns, in the worst case, as being proportional to
$m$ --- the number of unknowns, and to $n$ --- the number of
experimental points: ${\cal C}_{t} = 2Kmn$; effectively ${\cal C}_{t}
\sim {\cal O}(n)$, since $m=2$ is fixed. The factor $2$ comes from the
fact, that slicing is always two-sided. $K$, the proportionality
constant, is roughly equal to the number of bits in mantissa plus twice
the largest exponent used in floating point representation of real
numbers. One must remember, however, that a~single iteration, involving
all unknowns in turn, only rarely will suffice.  Fortunately, in the
linear case considered here, all necessary intervals, even if
calculated in a natural way (``naive''), have sharp ends, see Hansen's
\cite{Hansen} theorems on sharpness.  That is why the computed interval
enclosures are tight, also in the case of multivariate linear or
linearized models. The above statement is not necessarily true, when
the model is nonlinear.

\medskip\noindent
The spatial complexity of box slicing algorithm is very attractive.  At
any stage of the calculations we are working with at most $3$ interval
boxes (original, slice and rest), each of size proportional to the
number of unknowns $(m=2=const)$, so ${\cal C}_{s}\sim {\cal O}(1)$.

\subsection{Probing}

The purpose of probing is to determine whether the given box contains
the points with required properties, in our case -- the solutions.  The
boxes, which certainly do not contain at least one interesting point,
are eliminated from further considerations.

\begin{figure}[h]
 \epsfysize=0.7\hsize
 \epsfxsize=0.95\hsize
 \centerline{\epsfbox{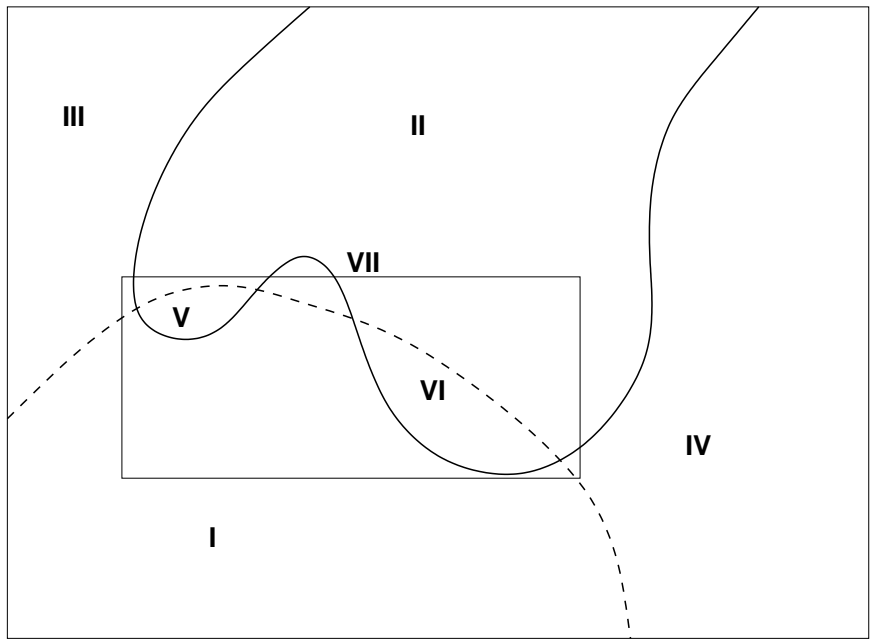}}
\caption{\sl The bounding box denotes initial search domain in the
plane $ab$.  The target, slightly oversized, is a rectangle located
near the center of figure.  Two remaining lines, solid and dashed,
divide the initial box into parts, for which some conditions are, or
are not, met.  The region bounded by dashed line is labelled as {\sf
I}, the region bounded by solid line has the label {\sf II}.  The
details concerning the remaining domains {\sf III} --- {\sf VII} are
given in text.
}
\label{inner_c}
\end{figure}

\medskip\noindent
Proving the existence of solutions within the probed domain is,
generally, not straightforward. Therefore, during probing, we will
rather seek every opportunity to discard the (sub)box under study. The
tests (``questions'') we are going to use during probing have to be
carefully selected, since
\begin{quotation}
\ldots\  {\sl probing has a logic of its own.}
\end{quotation}
\rightline{\sl M.H. van~Emden in \cite{Emden}}

\bigskip\noindent
Suppose that $p<q$ and we have obtained for some system of inequalities
${\cal I}$:
\begin{itemize}
\item
${\cal I}$ is non-failed for $x\le q\quad$ and
\item
${\cal I}$ is non-failed for $x\ge p$.
\end{itemize}
Can we say anything about the localization of the solutions of ${\cal
I}$, especially within the interval $\left[p,\ q\right]$? No, but on
the other hand the result:
\begin{center}
	$\left({\cal I}\right.$ is failed for $\left.x\le
	q\right)\quad$ and $\quad\left({\cal I}\right.$ is failed for
	$\left.x\ge p\right)$
\end{center}
\noindent
is a proof, that ${\cal I}$ has no solutions at all, while
\begin{center}
$\left({\cal I}\right.$ is failed for $\left. x\ge q\right)\quad $ and
$\quad \left({\cal I}\right.$ is failed for $\left.x\le p\right)$
\end{center}
implies, that the solutions, if any, must be located within the
interval $\left[p,\ q\right]$ (there are no solutions outside this
interval).

\medskip\noindent
Consider the Fig.\ref{inner_c}.  Think of the situation shown there as
applicable to just one experimental point, say first, and to all the
inequalities, in which this point is explicitly involved.  For sake of
simplicity from now on we drop the index numbering experimental points. 
Let's concentrate first on the united solution set.  We are going to
find such ${\mathbf a}$ and ${\mathbf b}$ that
\begin{equation}\label{united_cond}
	\left({\mathbf a}{\mathbf x} + {\mathbf b}\right) \cap {\mathbf
	y} \ne \emptyset
\end{equation}
for any pair $\left(x,\ y\right)\in \left({\mathbf x},\ {\mathbf
y}\right)$.  The opposite condition, when intervals ${\mathbf
a}{\mathbf x} + {\mathbf b}$ and ${\mathbf y}$ {\em are\/} disjoint, is
of better value for our purposes.  It may be written in conventional
interval notation as
\begin{equation}\label{just_two}
	\bigl( \overline{{\mathbf a}{\mathbf x} + {\mathbf b}}\
	<\ 
	\underline{\mathbf y} \bigr) \lor\
	\bigl( \underline{{\mathbf a}{\mathbf x} + {\mathbf b}}\
	>\
	\overline{\mathbf y} \bigr)
\end{equation}
The alternative (\ref{just_two}), written above,  provides us with the
correct answer to the question whether the condition
(\ref{united_cond}) is failed.  Boxes $\left({\mathbf a},\ {\mathbf
b}\right)$, for which (\ref{just_two}) is true, can be safely discarded
from further considerations, since all points belonging to them violate
(\ref{united_cond}).  All such boxes are localized outside the regions
marked as {\sf I} and {\sf II} in Fig.\ref{inner_c}, i.e. they may be
found in regions {\sf III}, {\sf IV} or {\sf VII}.  So, the alternative
(\ref{just_two}) may be used as a rejection criterion by box slicing
algorithm.  Think, however, what will happen, if the united set of
solutions is empty?  In such case we will finish with very small box,
still not sure whether or not there are any solutions in it.  In
contrast to ordinary point calculations, we cannot check it with a
single calculation --- there is still uncountable number of points
belonging to such a box.

\medskip\noindent
The solution of this dilemma is quite simple.  Starting with initial
box $V$ we discard those its parts, which satisfy the first term of the
alternative (\ref{just_two}), obtaining in result the box
$V_{d}\subseteq\ V$.  The box $V_{d}$ covers the region {\sf I} in
Fig.\ref{inner_c}.  Then we repeat the procedure, again starting with
$V$, but this time using the second term of (\ref{just_two}) as a
rejection criterion.  Now the resulting box is $V_{u} \subseteq V$,
covering the domain {\sf II} in Fig.\ref{inner_c}.  The solutions, if
any, must be located in the intersection $V_{d} \cap\ V_{u}$.  This
intersection, if not empty, becomes the new initial box for another
iteration.  Continuing this process we obtain better and better
interval hulls of regions marked as {\sf V} and {\sf VI} in
Fig.\ref{inner_c}. Procedure terminates (``eventually stabilizes'' in
language of \cite{Apt}), when $V_{u} = V_{d} = V$, or, in other words,
when the slicing becomes idempotent operation.  Empty intersection of
$V_{u}$ and $V_{d}$, at any stage of calculations, constitutes a proof
that the set of solutions is empty.

\medskip\noindent
So, asking the proper questions during probing is not a trivial thing. 
The main difficulty lies in construction of appropriate rejection
tests.  Only the tests $Q$ with property ($V \subset W \in $ {\mm
IR}$^{n}$):
\begin{equation}
	\biggl(Q(V)\ \hbox{\rm is non-failed}\biggr) \Longrightarrow\
	\biggl(Q(W)\ \hbox{\rm is non-failed}\biggr)
\end{equation}
are appropriate.  This is because we usually start with severely
overestimated initial box known to contain solutions.  We don't want to
have it rejected as a whole, in result of the first applied test.

\subsection{Rejection tests and search regions for various kinds of
solutions}
\subsubsection{United solutions}

In case of just two different measurements ($m=2$, ${\mathbf
x}_{1}\cap{\mathbf x}_{2}=\emptyset$) the problem may be quickly solved
``analytically'':
\begin{eqnarray}
	{\mathbf a} &=& \frac{{\mathbf y}_{2} - {\mathbf y}_{1}}{{\mathbf
	x}_{2} - {\mathbf x}_{1}} \label{simple_2pa}\\
	{\mathbf b} &=& \bigl({\mathbf y}_{1} - {\mathbf a}{\mathbf
	x}_{1}\bigr) \cap \bigl({\mathbf y}_{2} - {\mathbf a}{\mathbf
	x}_{2}\bigr)\label{simple_2pb}
\end{eqnarray}
Author cannot resist the temptation to comment on the elegance of the
expression (\ref{simple_2pb}) -- both measurements are treated in a
perfectly symmetric manner.  Despite this elegance, the expression for
${\mathbf b}$ does not necessarily produce the tight interval enclosure
for $b$, while the enclosure (\ref{simple_2pa}) for $a$ {\em is\/}
tight (sharp in Hansen's terminology, see \cite{Hansen}).  This is due
to the well known dependency problem, which is unavoidable here; namely
{\bf b} is expressed by, among other, ${\mathbf a}$, ${\mathbf x}_{1}$
and ${\mathbf y}_{1}$, while ${\mathbf a}$ was already calculated using
the same variables.  The lack of tightness is not a real problem here,
we can rectify it using the box slicing algorithm, which, by its
construction, always produces tight interval enclosures in linear case.

\medskip\noindent
The construction given by formulae (\ref{simple_2pa}) and
(\ref{simple_2pb}) can be used for establishing the bounds for the
initial box $V$, in which we will search for solutions (of any kind) of
the original linear problem.  It is sufficient to set
\begin{equation}\label{V_hull}
	V = \left({\mathbf a},\ {\mathbf b}\right) =
	\underline{\bigcup}_{jk}\ \left({\mathbf a}_{jk},\ 
	{\mathbf b}_{jk}\right)
\end{equation}
where ${\mathbf a}_{jk}$ and ${\mathbf b}_{jk}$ are intervals obtained
using values from measurements $j$ and $k$.  The convex hull written
above covers results for all pairs $(j,k)$ of experimental data
satisfying ${\mathbf x}_{j}\cap{\mathbf x}_{k}=\emptyset$.  Calculation
of (\ref{V_hull}) has thus the complexity ${\cal O}\left(n^{2}\right)$,
what may seem excessive, especially when the number of measurements is
large.  Instead we may start with initial box defined as
$V=\bigl(\left[-\omega,\ +\omega\right],\ \left[-\omega,\
+\omega\right]\bigr)$, where $\omega$ is some sufficiently large
number, say $10^{40}$, at a price of increased number of iterations
later.

\medskip\noindent
In summary: in order to find {\bf united solutions} we need the initial
box given by (\ref{V_hull}) and a pair of rejection rules given in
(\ref{just_two}).

\subsubsection{Tolerable solutions}
Those, if exist, are a subset of united solutions.  For this reason the
initial box might be constructed the same way as for united solutions,
see (\ref{V_hull}).  What we need are the rejection rules appropriate
for this case.  Rules (\ref{just_two}) are not sufficient.  We need to
discard boxes having no common parts with any tolerable solution, not
just with the experimental points.  Then we have to find at least one
tolerable solution.  Unfortunately, this is bad idea.  Driven by such a
condition, the procedure will converge not to the convex hull of
tolerable solutions, but to this single, specific solution instead. 
Besides, we still don't know how to find the first tolerable solution.

\medskip\noindent
Let us start with calculating
\begin{equation}\label{V_tol}
	V^{tol} = \left({\mathbf a}^{tol},\ {\mathbf b}^{tol}\right) =
	\left(\bigcap_{jk}\ {\mathbf a}_{jk},\ \bigcap_{jk}\ {\mathbf
	b}_{jk}\right)
\end{equation}
what can be done in the same loop, in which (\ref{V_hull}) is
calculated.  The box $V^{tol}$, if not empty, contains all straight
lines having common points with each experimental rectangle.  If it is
empty, then the tolerable solutions cannot exist.  Let's also calculate
the auxiliary intervals, one for each measurement
\begin{equation}\label{prime_tol}
	{\mathbf y}^{\prime}_{k} = {\mathbf a}^{tol}{\mathbf x}_{k} +
	{\mathbf b}^{tol}
\end{equation}
Now the rejection rules may be expressed as follows. Reject box
$\left({\mathbf a},\ {\mathbf b}\right)$, if for at least one
measurement, say $k$,
\begin{equation}\label{Rej_tol_con}
	\left({\mathbf a}{\mathbf x}_{k} + {\mathbf b}_{k}\right) \cap
	{\mathbf y}^{\prime}_{k} = \emptyset
\end{equation}
When calculating $V_{d}$ reject also boxes satisfying the condition
\begin{equation}
	\bigl(\overline{{\mathbf a}\underline{{\mathbf x}_{k}} +
	{\mathbf b}}\ <\ \underline{{\mathbf y}_{k}}\bigr)\
	\lor\ \bigl(\overline{{\mathbf a}\overline{{\mathbf x}_{k}} +
	{\mathbf b}}\ <\ \underline{{\mathbf y}_{k}}\bigr)
\end{equation}
for at least one measurement, while for $V_{u}$ use as a rejection
criterion
\begin{equation}
	\bigl(\underline{{\mathbf a}\underline{{\mathbf x}_{k}} +
	{\mathbf b}}\ >\ \overline{{\mathbf y}_{k}}\bigr)\ \lor\
	\bigl(\underline{{\mathbf a}\overline{{\mathbf x}_{k}} +
	{\mathbf b}}\ >\ \overline{{\mathbf y}_{k}}\bigr)
\end{equation}

\subsubsection{Controllable solutions}
Controllable solutions, if exist, are subset of united solutions.  So,
the starting box can be constructed using (\ref{V_hull}) again. 
Following the prescription (\ref{V_tol}) we construct
\begin{equation}
	V^{con} = \left({\mathbf a}^{con},\ {\mathbf b}^{con}\right)
\end{equation}
and then the auxiliary intervals
\begin{equation}
	{\mathbf y}^{\prime}_{k} = {\mathbf a}^{con}{\mathbf x}_{k} +
	{\mathbf b}^{con}
\end{equation}
If $V^{con}=\emptyset$ then the controllable solutions cannot exists. 
We discard all boxes $\left({\mathbf a},\ {\mathbf b}\right)$
satisfying (\ref{Rej_tol_con}) for at least one measurement. 
Additionally, when dealing with $V_{d}$ we use
\begin{equation}
	\min \left( \underline{{\mathbf a}\underline{{\mathbf x}_{k}} +
	{\mathbf b}},\ \underline{{\mathbf a}\overline{{\mathbf x}_{k}}
	+ {\mathbf b}}\right)\ >\ \underline{{\mathbf y}_{k}}
\end{equation}
and when slicing $V_{u}$ the following condition
\begin{equation}
	\max \left( \overline{{\mathbf a}\underline{{\mathbf x}_{k}} +
	{\mathbf b}},\ \overline{{\mathbf a}\overline{{\mathbf x}_{k}}
	+ {\mathbf b}}\right)\ <\ \overline{{\mathbf y}_{k}}
\end{equation}
as a rejection criterion.  It happens quite often that initially
$V^{con}\ne\emptyset$, but the set of controllable solution is
empty anyway.

\section{An Example}
Using a FORTRAN program, described in \cite{Joli}, I have fitted an
artificial data set consisting of $10$ points, with uncertainties in
both variables.  As the uncertainties for each experimental
measurement, $\sigma_{x}$ and $\sigma_{y}$, I have used third parts of
radii of the corresponding intervals, i.e. $\frac{1}{6}$th of their
widths.  The data are listed in Table \ref{data_t} and the results in
Table \ref{Joli_table}.

\begin{table}[h]
\caption{\sl Data used in exemplary calculations. $\sigma_{x}$ and
$\sigma_{y}$ were taken as rounded third parts of the corresponding
radii.}
\label{data_t}
\begin{center}
\begin{tabular}{rrrrrr}
\hline
center of ${\mathbf x}$ & 
   radius of ${\mathbf x}$  &
                 $\sigma_{x}$  &
                      center of ${\mathbf y}$  &
                                        radius of ${\mathbf y}$ &
                                                         $\sigma_{y}$\\
\hline\hline
 0.9      & 0.1       & 0.333      &  3.65     & 0.45     &\phantom{xx} 0.150\\
 1.9      & 0.1       & 0.333      &  4.60     & 0.40     & 0.133\\
 2.9      & 0.1       & 0.333      &  5.65     & 0.22     & 0.073\\
 3.9      & 0.1       & 0.333      &  6.60     & 0.40     & 0.133\\
 5.4      & 0.1       & 0.333      &  8.00     & 0.50     & 0.167\\
 5.9      & 0.1       & 0.333      &  8.55     & 0.35     & 0.117\\
 6.9      & 0.1       & 0.333      &  9.60     & 0.50     & 0.167\\
 8.7      & 0.1       & 0.333      & 11.30     & 0.50     & 0.167\\
 9.1      & 0.1       & 0.333      & 12.75     & 0.55     & 0.183\\
10.1      & 0.1       & 0.333      & 13.70     & 0.30     & 01.00\\
\hline
\end{tabular}
\end{center}
\end{table}

\begin{table}[h]
\caption{\sl Results produced by program taken from
\protect{\cite{Joli}}, as recorded from computer screen, without
any rounding.
}
\label{Joli_table}
\begin{center}
\begin{tabular}{lr|lr}
\hline
parameter & value & parameter & value\\
\hline\hline
$a_{LSQ}$           & 1.08530271 & $\sigma_{a}$  & 0.0136506381 \\
$b_{LSQ}$           & 2.43730211 & $\sigma_{b}$  & 0.0823259652\\
\hline
\end{tabular}
\end{center}
\end{table}

\medskip\noindent
I have also obtained an interval hull of united solution set for this
problem, using full uncertainties for both variables.  The result is
(rounded outwards to five significant figures):
$$
\left({\mathbf a},\ {\mathbf b}\right) = \bigl(\left[ 1.02270,\
1.13159\right],\ \left[ 2.06840,\ 2.96827\right]\bigr)
$$
As expected, $a_{LSQ}\in{\mathbf a}$ and $b_{LSQ}\in{\mathbf b}$. 
Using those values, the Table \ref{results_t} was then prepared.  In
order to calculate the ``corridor of errors'' ($y_{LSQ}-3\sigma$ and
$y_{LSQ}+3\sigma$ in Table) for LSQ results, I have mechanically
adopted the widely used formula, taking as $\sigma$ the following quantity:
\begin{equation}
	\sigma = \sqrt{\left|\frac{\partial y}{\partial a}\sigma_{a}
	\right|^{2} + \left|\frac{\partial y}{\partial b}\sigma_{b}
	\right|^{2}} = \sqrt{\left|x\sigma_{a}\right|^{2} +
	\left|\sigma_{b}\right|^{2}}
\end{equation}



\begin{table}[h]
\caption{\sl Results of exemplary calculations.  The intervals were
rounded outwards to two significant figures, while the numbers labeled
as $y_{LSQ}\pm 3\sigma$, resulting from point calculations, were
rounded conventionally.  The last column contains the ratio of widths
of corresponding uncertainty estimates: interval to that produced by
Least SQuares method.  First and last row do not correspond to any
measurement, they are included to illustrate the extrapolation behavior
of both approaches.
}
\label{results_t}
\begin{center}{\small
\begin{tabular}{rrrrrrrr}
\hline
 $x$ & $\underline{y}$ &
             $\overline{y}$ &
                      $\underline{y}_{fit}$ &
                           $\overline{y}_{fit}$  & $y_{LSQ}-3\sigma$ &
                                                         $y_{LSQ}+3\sigma$ &
                                                   w$_{INT}$/w$_{LSQ}$\\
\hline\hline
 0.0 &       &       &  2.06 &  2.97 & 2.434 & 2.441 & 120\\
\hline
 0.8 &  3.20 &       &  2.88 &  3.88 & 3.105 & 3.506 & 2.494\\
 1.0 &       &  4.10 &  3.09 &  4.10 & 3.272 & 3.773 & 2.016\\
\hline
 1.8 &  4.20 &       &  3.90 &  5.01 & 3.943 & 4.839 & 1.239\\
 2.0 &       &  5.00 &  4.11 &  5.24 & 4.111 & 5.105 & 1.137\\
\hline
 2.8 &  5.43 &       &  4.93 &  6.14 & 4.781 & 6.171 & 0.871\\
 3.0 &       &  5.87 &  5.13 &  6.37 & 4.949 & 6.438 & 0.833\\
\hline
 3.8 &  6.20 &       &  5.95 &  7.27 & 5.620 & 7.503 & 0.701\\
 4.0 &       &  7.00 &  6.15 &  7.50 & 5.787 & 7.770 & 0.681\\
\hline
 5.3 &  7.50 &       &  7.48 &  8.97 & 6.877 & 9.502 & 0.568\\
 5.5 &       &  8.50 &  7.69 &  9.20 & 7.045 & 9.768 & 0.555\\
\hline
 5.8 &  8.20 &       &  8.00 &  9.54 & 7.296 & 10.168 & 0.536\\
 6.0 &       &  8.90 &  8.20 &  9.76 & 7.464 & 10.434 & 0.525\\
\hline
 6.8 &  9.10 &       &  9.02 & 10.67 & 8.135 & 11.500 & 0.490\\
 7.0 &       & 10.10 &  9.22 & 10.89 & 8.302 & 11.767 & 0.482\\
\hline
 8.6 & 10.80 &       & 10.86 & 12.70 & 9.644 & 13.898 & 0.433\\
 8.8 &       & 11.80 & 11.06 & 12.93 & 9.811 & 14.165 & 0.429\\
\hline
 9.0 & 12.20 &       & 11.27 & 13.16 & 9.979 & 14.431 & 0.425\\
 9.2 &       & 13.30 & 11.47 & 13.38 & 10.147& 14.698 & 0.420\\
\hline
10.0 & 13.40 &       & 12.29 & 14.29 & 10.817& 15.763 & 0.404\\
10.2 &       & 14.00 & 12.50 & 14.52 & 10.985& 16.030 & 0.400\\
\hline
20.0 &       &       & 22.52 & 25.60 & 19.200& 29.086 & 0.312\\
\hline
\end{tabular}}
\end{center}
\end{table}

\medskip\noindent
As can be seen from the Table \ref{results_t}, within the range, where
the experimental points were taken, the results are comparable. 
However, if one examines the behavior of both methods outside this
domain, then the predictive power of the interval method is clearly
superior.  The LSQ uncertainty estimate for $y$ at $x=0$ seems at least
unrealistic, while for $x=20$ it is more than three times wider
comparing with the extrapolation based on interval result.

\medskip\noindent
It should be noted that the problem of interpolation of interval valued
experimental data has already been investigated in \cite{Chenyi}
(Lagrange interpolating polynomials) and in \cite{Markov} (essentially
the decomposition into set of basis functions --- generalized
polynomials).  In both cases the uncertainty appeared only in one
(dependent) variable.  No conclusions concerning physical meaning or
validity of obtained parameters could be drawn.

\section{Other applications in experimental sciences}

\subsection{Detection of outliers}
Suppose, that we are trying to find the convex hull of united solutions
for some linear problem and it appears empty.  This means, that one or
more data points are {\em outliers\/}.  There are several methods,
heuristic as well as based on probability calculus, which make possible
the identification of outliers.  Note, that contrary to the interval
methods, the LSQ method always produces some estimates of unknown
parameters, regardless of the presence of outliers.  Outliers can be
most easily spotted when the data and the best fitted line are plotted
simultaneously on the same graph.  This is rather tedious task, if
performed manually.  In control applications, either industrial or in
laboratory, when the environment is noisy, the misreadings or data
transmission errors may occur quite frequently and go undetected.  This
may seriously affect the quality and reliability of the control
procedures, sometimes leading to disastreous effects or even
fatalities.

\medskip\noindent
The advantage of interval methods over traditional LSQ method is then
evident: no outlier can go unnoticed.  Nevertheless, the question of
its identification still remains.  We should mention here, that other
methods usually fail, when the data set contains more than a single
outlier, or there are just two outliers located one next to the other.

\medskip\noindent
Proposed is  the following procedure for outlier detection:  repeat
finding the united solution set for the experimental data using relaxed
conditions (\ref{united}), i.e. with dropped requirement, that the
condition is true for every measurement.  In other words, we will
search for {\bf crude\/} solutions mentioned earlier. In order to find
this set we will discard boxes $\left({\mathbf a},\ {\mathbf
b}\right)$, which fail condition (\ref{united}) for at least $k$
measurements, where index $k$, initially set to zero, numbers the
consecutive trials.  The value of index $k$, at which we succeed for
the first time, tells us whether there are any outliers in the
investigated data set and, eventually, how many of them are there. 
Their identification is immediate: they all do not satisfy the
condition (\ref{united}), evaluated with most recently obtained
intervals ${\mathbf a}$ and ${\mathbf b}$.  Note, that there is no need
for a priori knowledge, {\em which\/} measurements are suspected.

\medskip\noindent
The method, outlined above, should be very robust -- it should be able
to detect and identify quite large number of outliers, even if they
constitute up to $50$\% of all measurements.  I have performed only a
very limited number of tests, with only one or two outliers.  The tests
show, that the method is working as expected.  It is obvious, that it
always terminates, when only two measurements remain, at the latest.

\subsection{Finding the asymptotic straight line}

Sometimes we need to find, based on experimental information, the
equation of a line, which is asymptotic or tangent to the investigated
curve.  The interval methods, described in this paper, may be of some
help in eliminating the subjectivity of the person processing this sort
of experimental data.  Nowhere in the literature I was able to find a
prescription on how to deal with problems of this kind.  Yet they are
important and so is the reliable determination of their relevant
parameters, including the uncertainties.  Consider, for example, the
process described by formula:
\begin{equation}
	y(t) = A_{0} + \sum_{j=1}^{k} A_{j} \exp
	\left(-\frac{t}{\tau_{j}}\right)
\end{equation}
where the number of various subprocesses, $k$,  need not to be known
precisely in advance and the unknown relaxation times $\tau$'s are well
separated and ordered in increasing order.  Recording such a process,
especially near $t=0$, we obtain
\begin{equation}
	y(t) \approx y(0) - \frac{A_{1}}{\tau_{1}}t
\end{equation}
that is, the equation of a straight line, providing that the entire
process is initially dominated by subprocess with shortest relaxation
time.  This is an equation commonly considered when studying
(photo)chemical reaction's kinetics, radioactive decay, multilevel
relaxation and many other processes.  Similar expressions can be
obtained for the initial permeability of ferromagnetic materials or
susceptibility of paramagnetics.  In all such cases the interesting
physical parameters are hidden in the slope of the corresponding
straight line being tangent to the experimental curve.  How can we
obtain the reliable values of parameters for such a line?

\medskip\noindent
Suppose, that we have found the interval hull of united solutions
describing first $n$ experimental points $\left({\mathbf a}_{n},\
{\mathbf b}_{n}\right)$.  Observe now, what should happen when we
enrich the data set with one more measurement and try to find united
solutions again (and the experimental points follow the straight
line!).  It is obvious, that $\left({\mathbf a}_{n+1},\ {\mathbf
b}_{n+1}\right) \subseteq \left({\mathbf a}_{n},\ {\mathbf
b}_{n}\right)$, if only the measurements are correct. This simple
observation is a basis for the proposed method.  Starting with two
first measurements (in this case the united solutions always exist),
keep finding the united solutions, enlarging the set of measurements by
one at every next trial.  Finish when either $\left({\mathbf a}_{n+1},\
{\mathbf b}_{n+1}\right) \not\subseteq \left({\mathbf a}_{n},\ {\mathbf
b}_{n}\right)$ or $\left({\mathbf a}_{n+1},\ {\mathbf b}_{n+1}\right) =
\emptyset$, whatever happens first.  The solution is then
$\left({\mathbf a}_{n},\ {\mathbf b}_{n}\right)$.  Of course, the data
should be properly ordered first, either as increasing or decreasing
sequence of $x$'s.

\section{Summary}
I have demonstrated how the interval analysis can be helpful in
objective and reliable processing of experimental data.  I have
discussed various types of possible solutions of linear systems, which
may be useful and desired under different circumstances.  The methods
presented here easily handle classical cases, when the uncertainties
affect only the dependent variable, the values of independent variable
being exact, as well as those with uncertainties in both variables. 
There is no need for, always more or less arbitrary, weighting the
experimental data.

\section{Acknowledgments}
The approach presented here was been inspired in great extent by the
work of Krzysztof R.~Apt \cite{Apt}.  After reading his paper author
did realize, that the partial order existing in {\mm IR}, the one
generated by the inclusion, is the fundamental property of this set. 
Interval methods should depend on it, as is the case with the well
known interval Newton algorithm.  The interval part of presented
numerical results was obtained using {\sf INTLIB} \cite{Kearfott} --
free interval library written in {\sf FORTRAN~77}.

This work was done as part of author's statutory activities in
Institute of Physics, Polish Academy of Sciences.

\end{document}